\documentclass[12pt]{article}
\usepackage{amssymb,latexsym}
\usepackage{amsmath}
\def\RR{\mathbb{R}}

\begin{document}

%------------------------------------------------------------------------

\title{\bf A quadrature formula associated\\ with a univariate quadratic spline quasi-interpolant}

\author{ Paul Sablonni\`ere\\ INSA and IRMAR, Rennes}
\date{December 2005}
\maketitle

\begin{abstract}
We study a new simple quadrature rule based on integrating a $C^1$ quadratic spline
quasi-interpolant on a bounded interval. We give nodes and weights for uniform and non-uniform
partitions. We also give error estimates for smooth functions and we compare this formula with
Simpson's rule.
\end{abstract}
  
%********************************************INTRODUCTION******************************************

\section {Introduction}

We study a new simple quadrature formula (QF)  based on integrating a $C^1$ quadratic spline
quasi-interpolant (qi) on a bounded interval. It can be considered as the companion of Simpson's rule (as midpoint rule is the companion of trapezoidal rule) in the sense that the signs of errors for both QF are opposite in almost all examples.\\
Here is an outline of the paper: in section 2, we recall the definition and main properties of the quadratic spline quasi-interpolants which are used. In section 3, we give the weights of the QF. In sections 4 and 5, we give estimates of the quadrature error for smooth functions in the cases of non-uniform and uniform partitions. In the latter case, we develop a more deeper study by using the associated Peano kernel and we compare our results to those of Simpson's rule. Finally, we illustrate our results with two numerical examples.\\
Regarding some recent papers close to ours, e.g.  De Swardt \& De Villiers \cite{SV}, and Lampret \cite{Lam}, we use quadratic splines as the first one, however we do not consider interpolants, but quasi-interpolants (other QF based on spline quasi-interpolants are given in \cite{Sab3}) and our QF seems to have slightly better properties, at least than the Lacroix rule. The paper by Lampret is more oriented to the calculation of sums via the Euler-Maclaurin formula: it would be interesting to compare the performances our formula in this domain.

%**************************** UNIVARIATE QUADRATIC SPLINES & QIs ***************************

\section{Univariate quadratic splines and discrete\\ quasi-interpolants}

Let  $X=X_n=\{x_0, x_1, \ldots, x_n\}$ be a partition of a bounded interval $I=[a,b]$ , with $x_0=a$
and $x_n=b$. For $1\le i\le n$, let $h_i=x_i-x_{i-1}$ be the length of the subinterval
$I_i=[x_{i-1},x_i]$ and let $\Gamma=\Gamma_n=\{0,1,\ldots,n+1\}$. We denote by ${\cal S}_2(X)$ the
$n+2$-dimensional space of $C^1$  quadratic splines on the partition $X$.
A basis of this space is formed by the family of quadratic B-splines $\mathcal{B}=\{B_i,\; i\in\Gamma\}$, with triple knots
$a=x_{-2}=x_{-1}=x_0$ and $b=x_n=x_{n+1}=x_{n+2}$. With these notations, the support of $B_i$ is
$[x_{i-2},x_{i+1}]$ for all $i\in\Gamma$. Define the set of data sites (or Greville points):
$$
\Theta_n=\{\theta_i=\frac12(x_{i-1}+x_i)\;;  i\in\Gamma\}.
$$
Note that $\theta_0=x_0$ and $\theta_{n+1}=x_n$.
In \cite{Sab1}\cite{Sab2}, we have proved  the existence of a unique {\sl discrete quasi-interpolant} (abbrev. dQI)
of type
$$
Qf=\sum_{i=0}^{n+1} \mu_i(f) B_i
$$
whose discrete coefficient functionals are respectively $\mu_0(f)=f(x_0)$, 
$\mu_{n+1}(f)=f(x_n)$ and, for $1\le i\le n$ :
$$
\mu_i(f)=a_i f(\theta_{i-1})+b_i f(\theta_i)+c_i f(\theta_{i+1}),
$$
which is exact on the space $\Pi_2$ of quadratic polynomials. Using the B-spline expansion of
monomials and  setting  $e_k(x)=x^k$ for $k\ge 0$, we get  (see e.g. \cite{dB},\cite{Mar1},\cite{Sch}) :
$$
e_0=\sum_{i\in\Gamma} B_i,\;\; e_1=\sum_{i\in\Gamma}  \theta_i B_i,\;\; 
e_2=\sum_{i\in\Gamma} x_{i-1}x_i B_i,
$$
and writing the equations $Qe_k=e_k $ for $ 0\le k\le 2$, we obtain:
$$
a_i=-\frac{\sigma_i^2\sigma'_{i+1}}{\sigma_i+\sigma'_{i+1}} , \;\;b_i=1+\sigma_i\sigma'_{i+1},
\;\; c_i=-\frac{\sigma_i(\sigma'_{i+1})^2}{\sigma_i+\sigma'_{i+1}},
$$
with
$$
\sigma_i=\frac{h_i}{h_{i-1}+h_i},\quad \sigma'_i=1-\sigma_i=\frac{h_{i-1}}{h_{i-1}+h_i}.
$$
We can now define the {\sl fundamental functions} of ${\cal S}_2(X_n)$ associated with the
quasi-interpolant $Q$:
$$
\bar B_0=B_0+a_1B_1,\quad \bar B_{n+1}=c_nB_n+B_{n+1},
$$
and, for $1\le i\le n$:
$$
\tilde B_i=c_{i-1}B_{i-1}+b_iB_i+a_{i+1}B_{i+1}.
$$
They allow to express $Qf$ in the following shorter form:
$$
Qf=\sum_{i\in \Gamma}f(\theta_i)\bar B_i,
$$
and to show that the infinity norm of $Q$ is equal to the Chebyshev norm of its 
Lebesgue function:
$$
\Lambda_{Q}(x)=\sum_{i\in \Gamma}\vert \tilde B_i \vert.
$$
In [20], Marsden proved the existence of a unique {\sl Lagrange interpolant} $v\in {\cal S}_2(X)$
satisfying $v(\theta_i)=f(\theta_i)$ for $0\le i \le n+1$. He also proved the following interesting result:\\

{\bf Theorem 1.} The infinity norm of the Lagrange operator is uniformly
bounded by $2$, for any partition $X$ of $I$.\\

There exists a similar result for the dQI $Q$ defined above :\\

{\bf Theorem 2.}  For any partition $X$ of $I$, the infinity norm of $Q$ is uniformly bounded
by $3$.  If the partition is uniform,
one has $\Vert Q \Vert_{\infty}=\frac{305}{207}\approx 1.4734$.\\

{\bf proof : } For the uniform case, see e.g. \cite{Sab2}. The first part  of the proof is easy because we have
$$
\vert \bar B_i \vert\le (\vert c_{i-1}\vert+b_i+\vert a_{i+1}\vert)B_i,
$$
therefore, by summing on all indices, as B-splines sum to one :
$$
\sum \vert \bar B_i \vert\le \max(\vert  a_i \vert+ b_i +\vert c_i \vert)=\max(1+2(\vert a_i\vert+\vert c_i \vert))=\max(1+2\sigma_i\sigma'_i)\le 3,
$$
since $b_i=1-a_i-c_i=1+\vert a_i\vert+\vert c_i\vert$, and $\sigma_i, \sigma'_i\le 1$ $\square$\\

{\bf Remark.} The results of this section are also valid when $X$ contains some knots of
multiplicity 
$2$ or $3$. Assume first that $\xi=x_p=x_{p+1}$ is a {\sl double knot}, then $Qf$ is only continuous at that point.
Moreover as $h_{p+1}=0$, we have $supp(B_{p})=[x_{p-2},\xi]$, $supp(B_{p+1})=[\xi-h_{p},\xi+h_{p+2}]$, $supp(B_{p+2})=[\xi,x_{p+3}]$.
Similarly, as $\sigma_{p+1}=0$, we have $a_{p+1}=c_{p+1}=0$ and $b_{p+1}=1$, hence:
$$
Qf=\sum_{i=0}^{p}\mu_i(f)B_i+f(\xi)B_{p+1}+\sum_{i=p+2}^{n+1}\mu_i(f)B_i.
$$

Now, if $\eta=x_{q-1}=x_q=x_{q+1}$ is a {\sl triple knot}, then $Q_2f$ has a discontinuity at this point.
Assume that $f$ is itself discontinuous and admits left and right limits $f(\eta^-)$ and $f(\eta^+)$.
Then as $h_q=h_{q+1}=0$, we have $supp(B_{q})=[\eta-h_{q-1},\eta]$, with $B_q(\eta^-)=1$ and $B_q(\eta^+)=0$, while
$supp(B_{q+1})=[\eta,\eta+h_{q+2}]$, with $B_{q+1}(\eta^-)=0$ and $B_{q+1}(\eta^+)=1$. As $\sigma_{q}=\sigma_{q+1}=0$,
we get $a_q=a_{q+1}=c_q=c_{q+1}=0$ and $b_q=b_{q+1}=1$, hence:
$$
Qf=\sum_{i=0s}^{q-1}\mu_i(f)B_i+f(\eta^-)B_{q}+f(\eta^+)B_{q+1}+\sum_{i=q+2}^{n+1}\mu_i(f)B_i.
$$

Finally, from theorem 2 and standard arguments in approximation theory (see \cite{DVL})
we deduce :\\

{\bf Theorem 3.} There exists a constant $0<C< 1$ such that for all 
$f\in W^{3,\infty}(I)$ and for all partitions of $I$, with $h=\max h_i$,
$$
\Vert f-Qf \Vert_{\infty}\le C h^3 \Vert D^3f \Vert_{\infty}. 
$$

%***********************************  QUADRATURE FORMULA ***************************************

\section{Quadrature formula associated with $Q$}

In this section, we compute the weights of the quadrature formula (abbr. QF) associated with the dQI $Q$, in 
the general case (non-uniform partition of the interval) and in the interesting particular case of a uniform partition.

%***********************************3.1  NON UNIFORM PARTITION ***************************************

\subsection{Case of a non uniform partition}

The QF associated with $Q$ is of course obtained by integrating $Qf$:
$$
{\cal I}_Q(f,I)=\int_a^b Qf=f(x_0)\int_a^b B_0+\sum_{i=1}^n \mu_i(f) \int_a^b B_i+f(x_n) \int_a^b
B_{n+1}.
$$
We know that $w_i=\int_a^b B_i=\frac16 (h_{i-1}+4h_i+h_{i+1})$ for $2\le i\le n-1$ and $$
w_0=\frac13 h_1,\;w_1=\frac16 (3h_1+h_2),\; w_n=\frac16 (h_{n-1}+3h_n),\;w_{n+1}=\frac13 h_n.
$$
Finally, we get the following QF:
$$
{\cal I}_Q(f,I)=\sum_{i\in\Gamma}w_i\mu_i(f),
$$
which can also be expressed in the classical form:
$$
{\cal I}_Q(f,I)=\sum_{i\in\Gamma}\bar w_i f(\theta_i),
$$
with $\bar w_0=w_0+a_1w_1, \bar w_{n+1}=c_n w_n+w_{n+1}$, and for $1\le i\le n$ :
$$
\bar w_i =\int_a^b\tilde B_i=c_{i-1}w_{i-1}+b_iw_i+a_{i+1}w_{i+1}.
$$
The explicit expression of these weights is rather complicated. We only prove  the following result :\\

{\bf Theorem 4.} \\
1) For all partitions, the weights $\bar  w_i$ satisfy $\sum_{i\in\Gamma}\vert\bar  w_i\vert\le 3(b-a)$.\\
2) More specifically, if there exists a $r>0$ such that $\displaystyle \frac{1}{r}\le \frac{h_{i+1}}{h_i}\le r$ for all $1\le i\le n-1$, then we have the more precise upper bound :
$$
\sum_{i\in\Gamma}\vert\bar  w_i\vert\le  (b-a)(1+2(\frac{r}{r+1})^2)
$$
{\bf proof : }   $\sum_{i\in\Gamma}\vert\bar  w_i\vert=\sum_{i\in\Gamma}\vert \int \bar   B_i \vert\le
\int\sum_{i\in\Gamma}\vert \bar   B_i \vert\le \max(\vert  a_i \vert+ b_i +\vert c_i \vert)\int \sum_{i\in\Gamma}B_i\le(b-a)\max(\vert  a_i \vert+ b_i +\vert c_i \vert)\le 3(b-a)$, the last inequality being a consequence of theorem 2. In addition,  when $\frac{1}{r}\le \frac{h_{i+1}}{h_i}\le r$ for all $i$, then we immediately deduce 
that $\frac{1}{r+1}\le \sigma_i, \sigma'_i \le \frac{r}{r+1}$ and  $\vert  a_i \vert+ b_i +\vert c_i \vert=1+2(\vert  a_i \vert+\vert c_i \vert)=1+2\sigma_i\sigma'_i\le 1+2(\frac{r}{r+1})^2 $, whence the result
$\square$\\

{\bf Remark 1.} As $Q$ is exact on $\Pi_2$, we have $E_Q(f,I)={\cal I}(f,I)-{\cal I}_Q(f,I)=0$ for all $f\in
\Pi_2$, therefore we can deduce that  $E_Q(f,I)=O(h^3)$, with $h=\max h_i$ for smooth functions
$f$ (see section 5 below, theorem 5). However, when  the partition $X$ satisfy the {\sl symmetry property}
$x_i=x_{n-i}, \,0\le i\le n$ (e.g. Chebyshev knots or, more generally, 
zeros of classical orthogonal polynomials in $I=[-1,1]$, see e.g.\cite{Dav}\cite{Gau1}\cite{Gau2}\cite{KU}), then the weights are also symmetric and we also obtain $E_Q(e_3,I)=0$. In that case, the error is a $O(h^4)$ for $f\in W^{4,\infty} (I)$, where $h=\max h_i$.\\

{\bf Remark 2.}  In general, the weights $\bar w_i $ are positive. However, they can be negative when the partition is strongly non-uniform. Take for example the following sequence of knots :
$$
X=X_7=\{-1, -0.9, -0.3, -0.2, 0.5, 0.6, 0.95, 1\}
$$
whose associated sequence of steplengths is :
$$
h_X=\{0.1, 0.6, 0.1, 0.7, 0.1, 0.35, 0.05\}
$$
Here we have  $r=7$, i.e. $\displaystyle \frac{1}{7}\le \frac{h_{i+1}}{h_i}\le 7$ for all $i$, and theorem 4 gives the upper bound $\sum_{i\in\Gamma}\vert\bar  w_i\vert\le \frac{81}{16}\approx 5.06$. The sequence of weights is the following :
$$
\bar w=\{.0146, .0122, .7463, .-.0622, .8780, -.0257, .4287, .0007, .0074\}
$$
We observe that the weights $\bar w_3$ and $\bar w_5$ are negative : however, their absolute values are small. Moreover, we have $\sum_{i\in\Gamma}\vert\bar  w_i\vert=2.17<5.06$.
%***********************************3.2   UNIFORM PARTITION ***************************************

\subsection{Case of a uniform partition}

Assume that $n\ge 5$ (in order to avoid boundary effects) and that the partition $X_n$ is uniform:
$h_i=h=\frac{b-a}{n}$ for $1\le i\le n$. In that case, the weights can be computed  once for all
and we obtain, by setting $f_j=f(\theta_j), j\in
\Gamma$ (see e.g. \cite{Sab1}\cite{Sab2}\cite{Sab3}):
$$
{\cal I}_Q(f,I)=h\,[\,\frac19(f_0+f_{n+1})+\frac{7}{8}(f_1+f_n)+\frac{73}{72}(f_2+f_{n-1})+
\sum_{i=3}^{n-2}f_i\,].
$$
As for some Newton-Cotes QF, this QF has a better degree of precision:\\

{\bf Theorem 4.} The QF associated with the dQI $Q$ on a uniform partition of $I$ is exact
on the space $\Pi_3$ of cubic polynomials. Therefore $E(f,I)=O(h^4)$ for $f$ smooth enough.\\

{\bf proof : } This is due to the symmetry of weights and nodes with respect to the midpoint of $I$.
$\square$\\

%**************************4.  ERROR  ESTIMATES/NON-UNIFORM PARTITION *******************

\section{Error estimate  on a non-uniform partition}

As we have already seen in section 4, $E_Q(f,I)={\cal I}(f,I)-{\cal I}_Q(f,I)=O(h^3)$  when the
partition $X$ of $I=[a,b]$ is non-uniform, with $h=\max_{1\le i\le n} h_i$.
In this section, we give more specific results when $f\in W^{3,\infty}(I)$.
From theorem 3, we immediately deduce:\\

{\bf Theorem 5.} There exists a constant $0<C_3<b-a$ such that for all $f\in W^{3,\infty}(I)$ and
for all partitions $X$  of $I$, with $h=\max_{1\le i\le n} h_i$:
$$
\vert E_Q(f,I) \vert\le C_3 h^3 \Vert D^3f \Vert_{\infty}.
$$
{\bf proof : } we know (theorem 3) that there exists a constant $0< C< 1$ such that $\Vert f-Qf \Vert_{\infty}\le C h^3 \Vert D^3f
\Vert_{\infty}$. As $E_Q(f,I)={\cal I}(f-Qf)$, we can write:
$$
\vert E_Q(f,I) \vert\le \int_a^b \vert f-Qf \vert\le (b-a)\, C h^3  \Vert D^3f \Vert_{\infty}.
$$
Therefore we get the result with $C_3=(b-a) C$.
$\square$
%***********************************5. ERROR/UNIFORM PARTITION ***********************************

\section{Error estimate on a uniform partition}

Assume that the partition $X_n$ is uniform ($h_i=h=\frac{b-a}{n}$ for $1\le i\le n$) 
and that $f\in W^{4,\infty}(I)$. As ${\cal I}_Q(f,I)$ is exact on $\Pi_3$, the Peano kernel
theorem (se e.g. \cite{Dav}, chapter III) gives:\\

$$
E_Q(f,I)=\frac16 \int_a^b K(t) D^4 f(t) dt,
$$
where
$$
K(t)=\int_a^b [(x-t)_+^3-Q(.-t)_+^3] dx,
$$
and
$$
Q(.-t)_+^3=\sum_{i\in\Gamma} (\theta_i-t)_+^3 \bar   B_i.
$$
As $\int_a^b (x-t)_+^3 dx=\int_t^b (x-t)^3 dx=\frac14 (b-t)^4$, we obtain:
$$
K(t)=\frac14 (b-t)^4-\sum_{i\in \Gamma}\tilde w_i (\theta_i-t)_{+}^3.
$$
We see immediately that $K(a)=K(b)=0$. First, it is clear that $K(b)=0$ since $(\theta_i-b)_{+}^3=0$
for all $i\in\Gamma$. Second, let $p(x)=(x-a)^3$, then 
${\cal I}(p)=\frac14 (b-a)^4={\cal I}_Q(p)=\sum_{i\in \Gamma}\tilde w_i p(\theta_i)$ since
$p\in\Pi_3$, therefore 
$K(a)=\frac14 (b-a)^4-\sum_{i\in \Gamma}\tilde w_i(\theta_i-a)^3=0$.

%***********************************  PEANO KERNEL ***********************************

\subsection{Sign structure of the Peano kernel}

For the sake of simplicity, we now assume that $I=[0,1]$, then $h=\frac1n$ and
$\theta_i=(i-\frac12)\,h$ for $1\le i\le n$, with  $\theta_0=0$ and $\theta_{n+1}=1$. Therefore we
have:
$$
K(t)=\frac14(1-t)^4-h\,[\frac78
\left(\frac{h}{2}-t\right)_+^3+\frac{73}{72}\left(\frac{3h}{2}-t\right)_+^3+
\frac{73}{72}\left((n-\frac32)h-t\right)_+^3
$$
$$
+\frac78 \left((n-\frac12)h-t\right)_+^3+\frac19
(1-s)_+^3]+\sum_{i=3}^{n-2}\left((i-\frac12)h-t\right)_+^3.
$$ 
It is not difficult to prove that $K(1-t)=K(t)$. Actually, setting
$$
K(t)=\frac14(1-t)^4-\sum_{i=0}^{n+1}w_i(\theta_i-t)_+^3 
$$
and using the symmetry of nodes and weights, we get
$$
K(1-t)=\frac14t^4-\sum_{i=0}^{n+1}w_i(t-\theta_i)_+^3.
$$
Now, we observe that $(\theta_i-t)_+^3-(t-\theta_i)_+^3=(\theta_i-t)^3$ and that the cubic
polynomial $p(s)=(s-t)^3$ is exactly integrated by the QF, hence:
$$
\int_0^1p(s)ds=\sum_{i=0}^{n+1}w_ip(\theta_i)\Longleftrightarrow
\frac14[(1-t)^4-t^4]-\sum_{i=0}^{n+1}w_i(\theta_i-t)^3=0.
$$
The above properties imply that 
$$
K(t)-K(1-t)=\frac14[(1-t)^4-t^4]-\sum_{i=0}^{n+1}w_i[(\theta_i-t)_+^3-(t-\theta_i)_+^3]=0.
$$
In the interval $I_1=[0,\theta_1]=[0,\frac{h}{2}]$, we have:
$$
K(t)=K_1(t)=\frac14 t^4-\frac{h}{9}t^3=\frac14t^3\left(t-\frac{4h}{9}\right),\;\;
K'_1(t)=t^2\left(t-\frac{h}{3}\right).
$$
Therefore, it is clear that $K_1$ has a minimum 
$K_1(\frac{h}{3})=-\frac{h^4}{972}$, that $K_1(\frac{4h}{9})=0$ and 
$ K_1(\theta_1)=K_1(\frac{h}{2})=\frac{h^4}{576}$.
Moreover $K'_1(\frac{h}{2})=\frac{h^3}{24}$.\\

In the interval $I_2=[\theta_1,\theta_2]=[\frac{h}{2},\frac{3h}{2}]$, we have:
$$
K(t)=K_2(t)=K_1(t)-\frac{7h}{8} \left(t-\frac{h}{2}\right)^3,\;\;
K'_2(t)=K'_1(t)-\frac{21h}{8}\left(t-\frac{h}{2}\right)^2.
$$
As $K'_1(\theta_2)=K'_1(\frac{3h}{2})=\frac{21h^3}{8}$, we see that $K'_2(\theta_2)=0$,
therefore we can factorize:
$$
K'_2(t)=\left(t-\frac{3h}{2}\right)\left(t^2-\frac{35h}{24}t+\frac{7h^2}{16}\right),
$$
and we deduce that $K'_2=0$ for $t=\frac{3h}{2}$ and $t=\bar t=\frac{35+\sqrt{217}}{48}h\approx
1.036h$. Therefore $K_2$ has a maximum at the latter point and
$$
K_2(\bar t)\approx 0.03 h^4\le \frac{h^4}{32} ,
$$
moreover, we have $ K_2(\theta_1)=K_1(\theta_1)=\frac{h^4}{576}$,
$ K_2(\theta_2)=\frac{h^4}{64}$, and we can factorize:
$$
K_2(t)-\frac{h^4}{64}=\frac18\left(t-\frac{3h}{2}\right)^2\left(t^2-\frac{17h}{18}t+\frac{h^2}{6}\right).
$$
In the interval $I_3=[\theta_2,\theta_3]=[\frac{3h}{2},\frac{5h}{2}]$, we have:
$$
K(t)=K_3(t)=K_2(t)-\frac{73h}{72} \left(t-\frac{3h}{2}\right)^3,\;\;
K'_3(t)=K'_2(t)-\frac{73h}{24}\left(t-\frac{3h}{2}\right)^2.
$$
As $K'_2(\theta_2)=0$, we see that $K'_3(\theta_2)=0$,
therefore we can factorize:
$$
K'_3(t)=\left(t-\frac{3h}{2}\right)\left(t^2-\frac{9h}{2}t+5h^2\right)=
\left(t-\frac{3h}{2}\right)(t-2h)\left(t-\frac{5h}{2}\right),
$$
and we deduce that $K_3$ has one maximum and two minima:
$$
K_3(2h)=\frac{h^4}{32},\;\;
K_3\left(\frac{3h}{2}\right)=K_3\left(\frac{5h}{2}\right)=\frac{h^4}{64}.
$$
Finally, we also obtain the factorization:
$$
K_3(t)-\frac{h^4}{64}=\frac14 \left(t-\frac{3h}{2}\right)^2\left(t-\frac{5h}{2}\right)^2.
$$
Now, in the intervals $I_i=[\theta_{i-1},\theta_i]=[(i-\frac32)h,(i-\frac12)h]$, for 
$4\le i\le n-2$, if we assume that
$$
K_i(t)-\frac{h^4}{64}=\frac14 \left(t-(i-\frac32)h\right)^2\left(t-(i-\frac12)h\right)^2,
$$
as we have, by definition
$$
K_{i+1}(t)=K_i(t)-\left(t-(i-\frac12)h\right)^3,
$$
we immediately oˆbtain
$$
K_{i+1}(t)-\frac{h^4}{64}=\frac14\left(t-(i-\frac12)h\right)^2
\left[\left(t-(i-\frac32)h\right)^2-4\left(t-(i-\frac12)h\right)\right]
$$
$$
=\frac14\left(t-(i-\frac12)h\right)^2\frac14\left(t-(i+\frac12)h\right)^2.
$$
Therefore, we deduce that, in the interval $I_i$, $K=K_i$ has one maximum and two minima 
$$
K_i((i-1)h)=\frac{h^4}{32},\;\;
K_i\left((i-\frac32)h\right)=K_i\left((i-\frac12)h\right)=\frac{h^4}{64}.
$$
In the three last subintervals $I_{n-1}=[1-\frac{5h}{2},1-\frac{3h}{2}],
I_{n}=[1-\frac{3h}{2},1-\frac{h}{2}]$ and
$I_{n+1}=[1-\frac{h}{2},1]$, the behaviour of
$K$ is symmetrical of that one in the three first subintervals. 
To sum up, we obtain \\

{\bf Theorem 6.} The Peano kernel $K$ is negative in the two small subintervals
$J_1=[0,\frac{4h}{9}]$ and $J_3=[1-\frac{4h}{9},1]$ and  positive in the subinterval
$J_2=[\frac{4h}{9}, 1-\frac{4h}{9}]$.

%*********************************** Error estimate for the  QF ***********************************

\subsection{Error estimate for the quadrature formula}

Let $f\in C^4(I)$ be a given function, then from 
$$
E_Q(f,I)=\frac16 \int_0^1 K(t) D^4 f(t) dt
$$
and theorem 6, we deduce:\\

{\bf Theorem 7.} For any function $f\in C^4(I)$, there exists two points $c,\bar c\in I$ such that 
$$
E_Q(f,I)= \frac{23}{5760}h^4D^4f(c)-\frac{1}{192}h^5 D^4f(\bar c).
$$

{\bf proof : }  From the  mean value theorem, we know that there exists\\
$c_1\in J_1, c_2\in J_2$ and $c_3\in J_3$ such that
$$
E_Q(f,I)=\frac16\left[D^4f(c_1)\int_0^{\frac{4h}{9}}K(t)dt+
D^4f(c_2)\int_{\frac{4h}{9}}^{1-\frac{4h}{9}}K(t)dt+D^4f(c_3)\int_{1-\frac{4h}{9}}^1 K(t)dt\right].
$$
We compute successively
$$
\int_0^{\frac{4h}{9}}K(t)dt=\frac14\int_0^{\frac{4h}{9}}t^3\left(t-\frac{4h}{9}\right)dt=
-\frac{64h^5}{295245},
$$
$$
\int_{\frac{4h}{9}}^{\frac{h}{2}}K(t)dt=\int_{\frac{4h}{9}}^{\frac{h}{2}}t^3\left(t-\frac{4h}{9}\right)dt=
\frac{1631h^5}{37791360},
$$
$$
\int_{\frac{h}{2}}^{\frac{3h}{2}}K(t)dt=\frac{59h^5}{2880},
$$
and for $2\le i\le n-2$:
$$
\int_{(i-\frac12)h}^{(i+\frac12)h} K(t)dt=\frac{23h^5}{960}.
$$
Then, using $nh=1$ and setting
$\gamma_1=\frac{64}{295245},\gamma_2=\frac{23}{960},\gamma_3=\frac{291149}{9447840}$,
we obtain successively:
$$
\int_{\frac32h}^{1-\frac52h} K(t)dt=(n-3)\frac{23h^5}{960}=\frac{23h^4}{960}-\frac{23h^5}{320}.
$$
$$
\int_{J_2}K(t)dt=\gamma_2h^4-h^5\left(2\times\frac{1631}{37791360}+
2\times\frac{59}{2880}-\frac{23}{320}\right)=\gamma_2h^4-\gamma_3h^5,
$$
$$
\int_{J_1}K(t)dt=\int_{J_3}K(t)dt=-\gamma_1h^5.
$$
Now, setting $\bar\gamma_j=\gamma_j/6$ for $1\le j\le 3$, we have:
$$
E_Q(f,I)=\frac16 \int_0^1 K(t) D^4 f(t) dt=\bar\gamma_2h^4 D^4f(c_2)-h^5
\left[\bar\gamma_1 D^4f(c_1)+\bar\gamma_3 D^4f(c_2)+\bar\gamma_1 D^4f(c_3)\right].
$$
Finally, setting $\bar \gamma_4=2\bar\gamma_1+\bar\gamma_3$,  the mean-value theorem implies that
there exists a point $c_4\in I$ such that 
$$
E_Q(f,I)=\bar\gamma_2 h^4D^4f(c_2)-\bar\gamma_4h^5 D^4f(c_4),
$$
with $\bar\gamma_2=\frac{23}{5760}$ and $\bar\gamma_4=\frac{1}{192}$.
$\square$\\

{\bf Remark.} For $n$ large enough, one can write
$$
E_Q(f,I)= \frac{23}{5760}h^4D^4f(c)+O(h^5),
$$
therefore, the main part of the error is contained in the first term.

%***********************************5.3  Comparison with Simpson ***********************************

\subsection{Comparison with composite Simpson's rule and extrapolation}

Let us compare the above error with the quadrature error $E_S(f,I)={\cal I}(f,I)-{\cal I}_S(f,I)$
of composite Simpson's rule. Taking $n=2m$ even, the latter is based on the $n+1$ points $\{x_j,
0\le j\le n\}$, while our QF is based on the set
$\{\theta_i, 0\le i\le n+1\}$ which has $n+2$ points. It is well known (see e.g. \cite{Dav}\cite{DR}\cite{Eng}\cite{KU}) that for any
function $f\in C^4(I)$, there exists a point $d\in I$ such that:
$$
E_S(f,I)=-\frac{1}{180}h^4D^4f(d).
$$
Therefore, {\sl when} $D^4f$ {\sl is of one sign} over the interval $I$, the quadrature errors are
of {\sl opposite signs}, i.e. the two formulas give {\sl lower and upper estimates} of the value of
the integral.\\

{\bf Remark:} as shown by numerical examples below, the following linear combination obtained by
extrapolation
$$
{\cal I}_{QS}(f,I)=\frac{1}{55}(32 {\cal I}_Q(f,I)+23 {\cal I}_S(f,I))
$$
gives a still better approximation, with an order $O(h^5)$ (notice the nice symmetry of weights).

%********************************6. NUMERICAL RESULTS **************************************

\section{Numerical results}

For the two following functions $f_1,f_2$ and $f_3$ (other examples are given in \cite{Sab3}),  we give the quadrature errors $E_Q(f)$, $E_S(f)$ and $E_{QS}(f)$ in terms of the number $n$ of subintervals:\\

%------------------------------------------------ Example 1 -----------------------------------
{\bf Example 1:} \\

${\cal I}(f_1,I)=\int_0^1 16 x^{3/2}sin\,(x^2)\,dx=3.2523064663781227544.$
$$
\begin{matrix}
n& E_Q(f_1)& E_S(f_1) & E_{QS}(f_1)\cr
 & & & \cr
64 & -.86(-7) & 1.23 (-7) &1.13 (-9)\cr
128 & -.54(-8) & .76 (-8) & .16 (-10)\cr
256 &  -.34(-9) &.47 (-9) & -.40 (-12)\cr
512 & -.21(-10) & .29 (-10) & -.52 (-13)\cr
1024 & -.13(-11) & .18 (-11) & -.33(-14)\cr
\end{matrix}
$$

%------------------------------------------------ Example 2 -----------------------------------

{\bf Example 2:} \\

$\displaystyle {\cal I}(f_2,I)=\int_0^1
\left( \frac{1}{(x-0.3)^2+0.01}+\frac{0.8}{(x-0.7)^2+0.04}\right)\,dx=35.880612010038328566$

\medskip

$$
\begin{matrix}
n& E_Q(f_2)& E_S(f_2) & E_{QS}(f_2)\cr
 & & & \cr
64 & -.19 (-5) & .23 (-5) & -.14 (-6)\cr
128 & -.11 (-6) & .14 (-6) & -.37 (-8)\cr
256 &  -.67 (-8) & .90 (-8) & -.11 (-9)\cr
512 & -.41 (-9) & .56 (-9) & -.35 (-11)\cr
1024 & -.25 (-10) & .35 (-10) & -.11 (-12)\cr
\end{matrix}
$$
\newpage
%------------------------------------------------ Example 3 -----------------------------------

{\bf Example 3 :} \\

$\displaystyle {\cal I}(f_3,I)=\int_{-1}^1
\frac{1}{1+16x^2}\,dx=0.6629088318340162325296195.$

\medskip

$$
\begin{matrix}
n& E_Q(f_2)& E_S(f_2) & E_{QS}(f_2)\cr
 & & & \cr
256 & -.33 (-10) & .46 (-10) & -.44 (-12)\cr
512 & -.21 (-11) & .28 (-11) & -.13 (-13)\cr
1024 &  -.13 (-12) & .18 (-12) & -.42 (-15)\cr
2048 & -.80 (-14) & .11 (-13) & -.13 (-16)\cr
4096 & -.50 (-15) & .69 (-15) & -.41 (-18)\cr
\end{matrix}
$$\\

Note that the signs of errors are opposite when $h$ is small enough.\\

%******************************** REFERENCES *******************************************

{\sl Abbreviations of editors in references:}

AP=Academic Press, London, New-York.
BV=Birkh\"auser Verlag, Basel.\\
Bl=Blaisdell, Waltham.
CUP=Cambridge University Press.\\
Dov=Dover Publ., New-York. JWS=John Wiley \& Sons, New-York.\\
IRMAR=Institut de Recherche Math\'ematique de Rennes.\\
SIAM=Society for Industrial and Applied Mathematics, Philadelphia.\\
SV=Springer-Verlag, Berlin, New-York.

%-----------------------------------------------------------------------------
\bigskip
{\bf AMS classification:} 41A15, 65D07, 65D25, 65D32.\\
{\bf Keywords : } Numerical integration, spline approximation.\\

{\bf Author's address: }

Paul Sablonni\`ere,\\
Centre de math\'ematiques, INSA de Rennes,\\
20 avenue des Buttes de Co\"esmes, CS 14315,\\
F-35043-RENNES C\'edex,\\
France\\

e-mail: psablonn@insa-rennes.fr

\end{document}